\documentclass[11pt,a4paper]{article}
\usepackage[T1]{fontenc}
\usepackage{lmodern,microtype}
\usepackage{amsmath,amssymb,amsthm,mathtools}
\usepackage{xcolor}
\usepackage{geometry,enumitem,needspace,hyperref}
\hypersetup{colorlinks=true,linkcolor=blue!50!black,citecolor=blue!50!black,urlcolor=blue!50!black,
 pdfauthor={NingNing Peng}}
\setlist{nosep,leftmargin=2.2em,topsep=4pt}
\theoremstyle{plain}
\newtheorem{theorem}{Theorem}[section]
\newtheorem{lemma}[theorem]{Lemma}
\newtheorem{proposition}[theorem]{Proposition}
\newtheorem{corollary}[theorem]{Corollary}
\theoremstyle{definition}
\newtheorem{definition}[theorem]{Definition}
\theoremstyle{remark}
\newtheorem{remark}[theorem]{Remark}

\newcommand{\RCA}{\mathsf{RCA}_0}
\newcommand{\WKL}{\mathsf{WKL}_0}
\newcommand{\ACA}{\mathsf{ACA}_0}
\newcommand{\CNE}{\mathsf{CNE}}

\newcommand{\Prob}{\mathcal P}
\newcommand{\Leb}{\ell}
\newcommand{\norm}[1]{\left\lVert#1\right\rVert}
\newcommand{\Nrm}{\mathcal N}
\newcommand{\ran}{\operatorname{ran}}

\newcommand{\eps}{\varepsilon}

\newcommand{\pair}[1]{\langle#1\rangle}
\numberwithin{equation}{section}

\author{
	Weiguang Peng, \\ School of Mathematics and Statistics, Southwest University, Chongqing, China
	\and
	NingNing Peng \\ School of Mathematics and Statistics, Wuhan University of Technology, Wuhan, China
}
\date{}
\newcommand{\Aclass}{\mathcal A}

\title{The Reverse Mathematics of Mixed Nash Equilibria in Continuous Games}
\hypersetup{pdftitle={The Reverse Mathematics of Mixed Nash Equilibria in Continuous Games}}
\begin{document}
\maketitle
\begin{abstract}
We study the existence of mixed Nash equilibria for continuous games in
second-order arithmetic. When compact presentations of the action spaces
are given jointly, we prove over $\RCA$ that equilibrium existence is
equivalent to $\WKL$. The reversal holds for one-player games on
Cantor space with $1$-Lipschitz payoffs. A separate corollary gives the
same equivalence for uniformly presented sequences of games. We also
obtain the reversal for symmetric two-player zero-sum games on the unit
interval whose payoffs are entire functions with rational Taylor
coefficients and fixed analytic bounds; in this case, a symmetric
equilibrium may be required. The compactness argument extends to games
with countably many players. We obtain rational finite-support
approximate equilibria in $\RCA$, including simultaneous approximations
for all players under a common continuity modulus. An $\omega$-model
consequence shows that equilibrium existence with separately asserted
compactness does not imply $\ACA$.
\end{abstract}
\noindent\textbf{Mathematics Subject Classification (2020).}
03B30, 03F35, 91A10, 91A06.
\par\noindent\textbf{Key words.} Reverse mathematics, continuous games,
mixed Nash equilibria, weak K\"onig's lemma, probability measures,
entire functions.

\section{Introduction}

Nash proved that every game with finitely many players and finitely
many pure strategies has a mixed Nash equilibrium \cite{nash50}.
Glicksberg extended mixed-equilibrium existence to games with finitely
many players, compact Hausdorff action spaces, and continuous payoff functions
\cite{glicksberg52}. In the compact metric case considered here, a
mixed strategy is represented by a probability measure on an action
space, and a mixed equilibrium is a profile of such measures under
which no player can increase the expected payoff by a unilateral
deviation. These existence theorems lead to a question in reverse
mathematics: which set-existence axioms are required to construct an
equilibrium from codes for the action spaces and payoffs?

Fixed point theorems provide a natural starting point. Shioji and
Tanaka \cite{st90} studied fixed point theory in weak subsystems of
second-order arithmetic. The usual reduction of finite Nash
equilibrium existence to Brouwer's fixed point theorem gives an upper
bound of $\WKL$ for finite games with explicitly given action lists;
see also \cite{ppy}. Peng and Yamazaki \cite{py2017} proved that
both the individual and sequential forms of Brouwer's theorem on
finite convex hulls are equivalent to $\WKL$ over $\RCA$.
They also proved in $\WKL$ the existence of pure Nash equilibria
for interval games whose payoffs are concave in each player's own
action. The latter result concerns a restricted class of games;
arbitrary continuous payoffs on compact action spaces require the
more general mixed-strategy formulation.


Peng, Peng and Yamazaki \cite{ppy} proved
that weak subsequence compactness of probability measures on compact
metric spaces is equivalent to $\ACA$. They used finite-game
approximations and a weakly convergent subsequence to establish mixed
equilibrium existence for continuous games in $\ACA$, and asked
whether equilibrium existence implies $\ACA$. Their argument
therefore leaves open whether the subsequence-compactness principle
is necessary for obtaining an equilibrium. The construction of a
single equilibrium and the extraction of a weakly convergent
subsequence from a prescribed sequence are different existence
problems, and their logical strengths must be examined separately.


Our first result determines the strength of equilibrium existence
under this joint-presentation convention. A game has a finite number
of players, with action spaces, payoff codes and continuity moduli
presented as in Definition~\ref{def:gamecode}. The number of players
is an input natural number. Let $\CNE$ denote equilibrium existence
for these games. The reversal needs only one player, a fixed action
space, and a fixed Lipschitz bound.

\Needspace{14\baselineskip}
\begin{theorem}[$\RCA$]\label{thm:main}
The following are equivalent.
\begin{enumerate}
\item $\WKL$.
\item $\CNE$.
\item Every one-player game on Cantor space $2^{\mathbb N}$ with a
$1$-Lipschitz payoff taking values in $[-1,0]$ has a mixed Nash
equilibrium.
\end{enumerate}
In the third assertion, Cantor space has the metric
\[
 d(x,y)=2^{-\min\{k\ge0:x(k)\ne y(k)\}}
 \quad(x\ne y),\qquad d(x,x)=0.
\]
\end{theorem}

The upper-bound proof treats all games in a uniformly presented
sequence in one compact coordinate space. This yields the following
sequential form, denoted by $\mathsf{SeqCNE}$.

\Needspace{7\baselineskip}
\begin{corollary}[Sequential existence; $\RCA$]\label{cor:sequential}
The following assertion is equivalent to $\WKL$: every uniformly
presented sequence of continuous games, with the data of
Definition~\ref{def:gamecode}, has a jointly coded sequence of mixed
Nash equilibria.
\end{corollary}

The joint output in Corollary~\ref{cor:sequential} is supplied by
the compactness argument itself. No countable choice principle is
used to assemble independently obtained equilibria. The assertion
concerns existence of a sequence and makes no claim of a computable
equilibrium-selection operation. When compactness is asserted
separately for each action space, the exact classification for an
internally variable finite number of players requires further
analysis. Nevertheless, Corollary~\ref{cor:separation} gives an
$\omega$-model in which equilibrium existence under that convention
holds and $\ACA$ fails. It therefore answers the implication
question of \cite{ppy} negatively without assuming a general
principle for collecting compactness witnesses.

Our second principal result retains the reversal under strong
regularity restrictions on a fixed interval. A two-player zero-sum
game with row payoff $A$ is symmetric when $A(x,y)=-A(y,x)$.
Let $\Aclass^{\mathbb Q}$ be the class of such entire payoffs with
exact rational Taylor coefficients and the common analytic bounds
in Definition~\ref{def:analyticclass}.

\Needspace{8\baselineskip}
\begin{theorem}[$\RCA$]\label{thm:skewequilibrium}
The following assertion is equivalent to $\WKL$: every zero-sum game
on $[0,1]^2$ with row payoff in $\Aclass^{\mathbb Q}$ has a symmetric
mixed Nash equilibrium $(\mu,\mu)$. The equivalence also holds if an
arbitrary mixed equilibrium is required, or if the rational-coefficient
restriction is omitted while the effective analytic representation
and bounds are retained.
\end{theorem}

For the upper bound in Theorem~\ref{thm:main}, we impose the
probability-functional identities and Nash inequalities as closed
conditions in a compact coordinate space. The finite Nash theorem
\cite{ppy} verifies finite satisfiability. This obtains an equilibrium
directly from compactness. The Cantor-space reversal extracts a path
through an infinite binary tree from a maximizing probability
functional. The proof of Theorem~\ref{thm:skewequilibrium} instead
uses an orthogonal family of rational polynomials to convert
equilibrium inequalities into equations for a signed tree flow.
Increasing orders of vanishing make each Taylor coefficient a finite
rational sum and preserve the prescribed analytic bounds.

The same compactness method applies to a single game with countably
many interacting players. We also obtain rational finite-support
approximate equilibria in $\RCA$. Individual continuity moduli
suffice for approximating the equilibrium conditions of any specified
finite set of players; a supplied common modulus yields simultaneous
approximation for all players. These assertions distinguish the
exact existence principle from its finite approximation consequences.

Section~\ref{sec:representations} gives the representations and explains
the role of joint compactness data. Sections~\ref{sec:compact}
and~\ref{sec:skewentire} prove the two theorems and the sequential
corollary. Section~\ref{sec:countably-many-players} treats countably
many players and approximate equilibria. The final section compares
equilibrium existence with weak compactness and states the remaining
questions. 
presentations may replace explicit finite-net sequences.

\section{Preliminaries}\label{sec:representations}

We use the subsystems of second-order arithmetic and the metric-space
representations of Simpson \cite{simpson}. Unless indicated otherwise,
all constructions are carried out in $\RCA$. Following
\cite{ppy}, a complete separable metric space $\widehat X$ is
represented by a countable dense metric space, and its points are fast
Cauchy sequences. A compact presentation includes a sequence of finite
$2^{-k}$-nets. The Banach space $C(\widehat X)$ is presented in the
uniform norm using the basic metric functions of \cite{ppy};
continuous-function codes include moduli of uniform continuity.

\begin{definition}\label{def:gamecode}
A continuous game consists of a natural number $n\ge1$, nonempty
compact complete separable metric action spaces $\widehat A_i$
$(i<n)$, and payoff functions
\[
 f_i\in C\left(\prod_{j<n}\widehat A_j\right)\qquad(i<n).
\]
The metric data, payoff codes and continuity moduli are given jointly.
The relevant unital Banach function-space presentations are also given
jointly: the unit and basic metric functions have specified names and
pointwise interpretations, and the supremum norms of their rational
linear combinations are uniformly available. These norm presentations
may be supplied for the factors or through the common tensor
presentation of $C(\prod_{j<n}\widehat A_j)$ in \cite{ppy}.
Equivalently, one may supply joint finite-net sequences for the action
spaces together with the payoff codes and moduli;  For a sequence of games, these data
are uniform in the game index.
\end{definition}

The number of players is an input natural number, rather than an
externally fixed integer. The significance of this convention is
discussed in \S\ref{sec:scope}.

A probability measure on $\widehat X$ is represented, as in
\cite{ppy}, by a positive normalized linear functional
\[
 \mu:C(\widehat X)\longrightarrow\mathbb R,
 \qquad \mu(1)=1,
 \qquad g\ge0\Longrightarrow\mu(g)\ge0.
\]
We also call $\mu$ a probability functional and write
$\Prob(\widehat X)$ for the space of such functionals. Equality is
extensional equality on continuous functions.

We use the effective Stone--Weierstrass correspondence and the
finite-product construction of Peng, Peng and Yamazaki \cite{ppy}
without reproving them. Both apply to the basic function presentations
in $\RCA$.
For a jointly given finite tuple $(\mu_i)_{i<n}$, the product is the
norm-one extension of
\[
 (\otimes_{i<n}\mu_i)
 \left(\sum_{r<m}\prod_{i<n}p_{ir}\right)
 =\sum_{r<m}\prod_{i<n}\mu_i(p_{ir}).
\]
In particular, evaluation on a finite tensor sum is polynomial in
finitely many functional coordinates.

For a mixed strategy profile $\mu=\pair{\mu_i:i<n}$, put
\[
 E_\mu(f_i)=(\mu_0\otimes\cdots\otimes\mu_{n-1})(f_i).
\]
Write $E_{a,\mu_{-i}}(f_i)$ for the expectation obtained by replacing
$\mu_i$ with the Dirac measure at $a$. Following
\cite{ppy}, the profile is a mixed Nash equilibrium if
\begin{equation}\label{eq:nash}
 E_\mu(f_i)\ge E_{a,\mu_{-i}}(f_i)
 \qquad(i<n,\ a\in\widehat A_i).
\end{equation}
Positivity gives the corresponding inequalities for mixed deviations.

\begin{definition}\label{def:principles}
The principle $\CNE$ asserts that every game in
Definition~\ref{def:gamecode} has a mixed Nash equilibrium.
The principle $\mathsf{SeqCNE}$ asserts that every uniformly presented
sequence of such games has a jointly coded sequence of mixed Nash
equilibria.
\end{definition}

No support code or assignment of real values to all Borel sets is
required. As observed in \cite{ppy}, for an open set $U$
the supremum
\[
 \sup\{\mu(g):0\le g\le1,\ g=0\text{ outside }U\}
\]
need not exist as a real number in $\RCA$. The functional
representation does not require such suprema.

\subsection{Uniformity of compactness data}\label{sec:scope}

Joint unital norm presentations
with specified basic metric functions yield joint finite-net sequences
in $\RCA$. Thus the two forms of compact presentation in
Definition~\ref{def:gamecode} are equivalent.

One can instead retain only metric and continuous-map codes and assert
compactness separately for each action space. For an internally
variable number $n$ of players, the assertion
\[
 \forall i<n\ \exists W_i\ \operatorname{Net}(A_i,W_i)
\]
does not by itself supply the joint witness required by
\[
 \exists W\ \forall i<n\ \operatorname{Net}(A_i,W_i),
\]
where $W_i$ is the $i$th section of $W$. For each externally fixed
finite number of players, existential elimination combines the
witnesses, so Theorem~\ref{thm:main} applies. Its reversal remains
valid after adjoining players with singleton action spaces.

A general collection principle for these witnesses is not provable in
$\WKL$. Applied to discrete spaces arising from a finite family of
nonempty bounded subsets of $\mathbb N$, joint $1/4$-nets give a common
bound for the family. This implies the finite-union-of-finite-sets
principle, which is equivalent over $\RCA$ to
$\mathsf B\Sigma^0_2$ and is not provable in $\WKL$
\cite{fm}. This rules out a general witness-collection argument; it does not
classify equilibrium existence under the separate-witness convention.

The question concerning $\ACA$ in \cite{ppy} nevertheless has a
negative answer under either convention.

\begin{corollary}\label{cor:separation}
There is an $\omega$-model of $\RCA$ in which every continuous game has a mixed Nash
equilibrium, even with separately asserted compactness, but $\ACA$
fails. Hence this equilibrium-existence assertion does not imply
$\ACA$ over $\RCA$.
\end{corollary}

\begin{proof}
The relativized low basis theorem yields a countable Scott ideal
$\mathcal I$ consisting of low sets such that
$M=(\mathbb N,\mathcal I)\models\WKL$ and $0'\notin\mathcal I$;
see \cite{js,simpson}. In this $\omega$-model, the number of players
in any game is externally finite. The compactness witnesses and payoff
codes therefore have a finite join in $\mathcal I$. When needed,
$\WKL$ supplies the uniform continuity moduli. Their finite join
gives the presentation of Definition~\ref{def:gamecode}, and
Theorem~\ref{thm:main}, applied in $M$, gives an equilibrium.
Since $0'\notin\mathcal I$, $M$ does not satisfy $\ACA$.
\end{proof}

We shall use the finite-intersection form of compactness in $\WKL$
\cite{simpson}, its equivalent separation principle for
disjoint $\Sigma^0_1$ sets \cite{simpson}, and Harrington's
$\Pi^1_1$-conservation theorem \cite{simpson}; see also
\cite{avigad}.

\section{Continuous games and weak K\"onig's lemma}\label{sec:compact}

\begin{lemma}[$\RCA$]\label{lem:coordinates}
From a compact presentation of $\widehat X$, one can uniformly obtain
a countable rational function algebra $P_X$ whose norm completion is
$C(\widehat X)$. For each $p\in P_X$, the uniform norm has a real
code and there is a rational bound $B_p\ge\max\{1,\norm p_\infty\}$.
The probability functionals are represented by effective closed
conditions in $\prod_{p\in P_X}[-B_p,B_p]$. These coordinates and the
original functional codes are uniformly interconvertible.
\end{lemma}

\begin{proof}
Take the rational algebra generated by the basic metric functions.
The effective correspondence of \cite{ppy} identifies its
norm completion with $C(\widehat X)$. Finite nets and continuity
moduli compute the uniform norms, while the algebraic expressions give
the bounds $B_p$.

On the coordinates $v(p)$ impose
\begin{equation}\label{eq:functionalconstraints}
 v(1)=1,\quad v(p+q)=v(p)+v(q),\quad
 v(cp)=cv(p),\quad |v(p)|\le\norm p_\infty
 \qquad(c\in\mathbb Q).
\end{equation}
These are effective closed conditions. They define a unital
contractive rational-linear functional. If $(p_m)$ is a fast uniform
approximation to $f$, set
\[
 \mu(f)=\lim_m v(p_m).
\]
The bound $|v(p_m)-v(p_l)|\le\norm{p_m-p_l}_\infty$ gives a prescribed
Cauchy rate and makes the extension independent of the name. It also
identifies distinct algebraic expressions of the same function.

The extension is positive. If $g\ge0$ and $M\ge\norm g_\infty$ is
rational, then $\norm{M-g}_\infty\le M$, whence
$M-\mu(g)=\mu(M-g)\le M$. Conversely, a positive normalized
functional is contractive and satisfies \eqref{eq:functionalconstraints}.
The constructions are uniform. Dividing the $p$th coordinate by $B_p$
places the ambient product in a compact Hilbert cube.
\end{proof}

\begin{proposition}[$\WKL$]\label{prop:upper}
Every uniformly presented sequence of continuous games has a sequence
of mixed Nash equilibria.
\end{proposition}

\begin{proof}
For one game, use Lemma~\ref{lem:coordinates} for all players.
By the effective approximation and product constructions of
\cite{ppy}, there are
uniformly given finite tensor sums $P_{i,k}$ with
$\norm{f_i-P_{i,k}}_\infty\le\eps_k=2^{-k}$. For each deviation $a$
in the specified dense subset of $\widehat A_i$, impose
\begin{equation}\label{eq:approxclosed}
 E_\mu(P_{i,k})\ge E_{a,\mu_{-i}}(P_{i,k})-2\eps_k.
\end{equation}
The product formula in \S\ref{sec:representations} makes this a
continuous polynomial inequality in finitely many coordinates.
Together with \eqref{eq:functionalconstraints}, these are effective
closed conditions in a compact Hilbert cube.

Consider a finite subfamily. For each player, take a finite list
containing all deviations mentioned in that subfamily, or one fixed
dense point if no deviation is mentioned. Use action labels, allowing
different labels to denote the same point. Restrict the original
payoffs to these lists. The finite Nash theorem in $\WKL$
\cite{ppy} supplies probability vectors
and hence atomic functionals
$\mu_i(g)=\sum_jw_{ij}g(a_{ij})$. They satisfy all functional
conditions. The finite-game inequalities and the two approximation
errors imply each selected instance of \eqref{eq:approxclosed}.

Compactness \cite{simpson} yields a common solution.
For every dense deviation and every $k$,
\[
 E_\mu(f_i)-E_{a,\mu_{-i}}(f_i)\ge-4\eps_k.
\]
Letting $k$ tend to infinity gives the Nash inequality, and continuity
extends it to all pure actions.

For a sequence of games, include the game index among the coordinates.
A finite subfamily of conditions involves finitely many games.
Combine their finite restrictions into one finite game with separate
blocks of players, each payoff depending only on its own block.
A finite Nash equilibrium verifies finite satisfiability. Compactness
in the single resulting coordinate space gives the sequence.
\end{proof}

\begin{proposition}[$\RCA$]\label{prop:lowercantor}
The third assertion of Theorem~\ref{thm:main} implies $\WKL$.
\end{proposition}

\begin{proof}
Let $T\subseteq2^{<\mathbb N}$ be an infinite prefix-closed tree.
On Cantor space define
\[
 h_T(x)=\sum_{m\ge1}2^{-m}\mathbf1[x\upharpoonright m\notin T],
 \qquad u_T(x)=-h_T(x).
\]
The tail after stage $N$ is at most $2^{-N}$. Each cylinder indicator
of length $m\ge1$ is a basic metric function, with centre
$\sigma0^\infty$, inner radius $2^{-m}$ and outer radius $2^{1-m}$.
Thus the partial sums give a uniform-norm code in the original
representation. If $x,y$ first differ at coordinate $k$, then
\[
 |u_T(x)-u_T(y)|\le\sum_{m>k}2^{-m}=2^{-k}=d(x,y).
\]
Hence $u_T$ is $1$-Lipschitz and takes values in $[-1,0]$.

For each $N$, extend a node of $T$ of length $N$ by zeros to obtain
$a$ with $h_T(a)\le2^{-N}$. If $\mu$ is a mixed equilibrium, then
$0\le\mu(h_T)\le2^{-N}$ for every $N$, so $\mu(h_T)=0$. For a
nonempty string $\sigma\notin T$,
\[
 2^{-|\sigma|}\mathbf1_{[\sigma]}\le h_T,
 \qquad \mu(\mathbf1_{[\sigma]})=0.
\]

The root cylinder has mass one. Starting with the lower bound
$\eta=1/2$, suppose that a cylinder has mass greater than $\eta$.
Approximate the masses of its two children within $\eta/8$.
At least one approximation exceeds $3\eta/8$; the corresponding
child has mass greater than $\eta/4$ and therefore belongs to $T$.
Choose that child and replace $\eta$ by $\eta/4$. Recursive
comprehension gives a path through $T$.
\end{proof}

\begin{proof}[Proof of Theorem~\ref{thm:main}]
Proposition~\ref{prop:upper}, applied to a constant sequence, gives
$(1)\Rightarrow(2)$. The canonical compact presentation of Cantor
space gives $(2)\Rightarrow(3)$. Proposition~\ref{prop:lowercantor}
gives $(3)\Rightarrow(1)$.
\end{proof}

\begin{proof}[Proof of Corollary~\ref{cor:sequential}]
Proposition~\ref{prop:upper} proves
$\WKL\Rightarrow\mathsf{SeqCNE}$.
Conversely, constant repetition of a game and extraction of the first
equilibrium give $\CNE$. Theorem~\ref{thm:main} then yields $\WKL$.
\end{proof}

\begin{remark}\label{rem:finite-input}
Peng, Peng and Yamazaki distinguish two representations of finite-game
sequences \cite{ppy}. Their $\ACA$
reversal uses $A_e=\{0\}\cup\{t+1:\alpha(t)=e\}$ for an injection
$\alpha$. Although each action set has at most two elements, a
sequence of complete lists would decide $\ran(\alpha)$. They also
observe the $\WKL$ upper bound when the requisite finite data are
given uniformly; see \cite{py2017}.
Proposition~\ref{prop:upper} extends this latter formulation to jointly
presented compact action spaces.
\end{remark}

\section{Symmetric equilibria for entire payoffs}\label{sec:skewentire}

Throughout this section both players have action space $[0,1]$, and the
payoffs are $A$ and $-A$. We write
$A(\mu,\nu)=(\mu\otimes\nu)(A)$, using the product functional of
\cite{ppy}, and put $\Leb(g)=\int_0^1g(x)\,dx$.
The latter functional is available in $\RCA$ for continuous functions
with supplied moduli \cite{simpson}.

\subsection{Entire-function codes}

For a rational polynomial $p=\sum_ka_kx^k$, put
\[
 |p|_{\mathrm c}=\sum_k|a_k|,
 \qquad \Nrm(p)=\frac{p}{1+|p|_{\mathrm c}}.
\]
Thus $|\Nrm(p)|_{\mathrm c}\le1$, and normalization preserves the
vanishing of moments. For a bivariate power series
$H(x,y)=\sum_{a,b\ge0}c_{ab}x^ay^b$, put
\[
 \norm H_R=\sum_{a,b\ge0}|c_{ab}|R^{a+b}.
\]
The corresponding notation is also used for univariate series.
An effective entire-function code consists of the coefficient array
and a joint modulus of absolute convergence, with respect to total
degree, at every positive integral radius $R$.

\begin{definition}\label{def:analyticclass}
Let $\Aclass$ be the class of effectively coded entire payoffs satisfying
\[
 A(x,y)=-A(y,x),\qquad
 \norm A_2\le2^{-10},\qquad
 \norm A_R\le2^{R^2+3}\quad(R\ge1\text{ an integer}),
\]
with the convergence moduli and these bounds included in the
presentation. Let $\Aclass^{\mathbb Q}$ be the subclass whose Taylor
coefficients are supplied as a joint array of exact rational numbers.
A symmetric mixed strategy profile has the form $(\mu,\mu)$.
\end{definition}

The skew-symmetry condition expresses symmetry of the zero-sum game,
since $-A(x,y)=A(y,x)$. The coefficient bounds give
$\norm{\partial_x^p\partial_y^qA}_\infty\le2^{-10}p!q!$ on $[0,1]^2$,
and hence common derivative and Lipschitz bounds for this class.

\begin{lemma}[$\RCA$]\label{lem:entire}
Suppose that rational polynomials $p_{j\ell},q_{j\ell}$, rational
numbers $\epsilon_{j\ell}$, and natural numbers $m_j,d_j$ are given
uniformly, with $m_j\le2$ and, for every $\ell<m_j$,
\[
 |p_{j\ell}|_{\mathrm c},|q_{j\ell}|_{\mathrm c}\le1,
 \qquad \deg p_{j\ell},\deg q_{j\ell}\le d_j,
 \qquad |\epsilon_{j\ell}|\le1.
\]
A zero polynomial may be assigned degree bound zero. Then
\[
 H(x,y)=\sum_{j\ge0}2^{-j-(d_j+4)^2}
       \sum_{\ell<m_j}\epsilon_{j\ell}p_{j\ell}(x)q_{j\ell}(y)
\]
has an effective entire-function code, uniformly in these data.
Writing $H_{<s}$ for the sum over $j<s$, for every integer $R\ge1$
we have
\begin{equation}\label{eq:entiretail}
 \norm{H-H_{<s}}_R\le2^{R^2-s+2},\qquad
 \norm{H-H_{<s}}_2\le2^{-s-14}.
\end{equation}
All derivatives have uniform codes and satisfy
\[
 \norm{\partial_x^p\partial_y^qH}_\infty\le2^{-14}p!q!
 \quad\text{on }[0,1]^2.
\]
\end{lemma}

\begin{proof}
A polynomial of degree at most $d$ and coefficient norm at most one
has norm at radius $R\ge1$ at most $R^d$. Each tensor summand therefore
has norm at most
\[
 2^{-j-(d+4)^2}R^{2d}
 \le2^{-j-(d+4)^2+2Rd}\le2^{R^2-j}.
\]
Summing over $j\ge s$ and at most two summands per stage proves the
first estimate. At $R=2$ the exponent is $-j-d^2-6d-16$, giving the
second estimate.

Each coefficient is the limit of rational partial sums with the stated
rate. Given a required Taylor-tail bound at radius $R$, choose $s\ge1$
so that the stage tail is sufficiently small and put
$M=2\max_{j<s}d_j$. Terms of total degree greater than $M$ receive
contributions only from stages $j\ge s$. This gives a joint modulus
for Taylor truncation.

For a finite polynomial $P$, coefficientwise differentiation gives
\[
 \norm{\partial_x^p\partial_y^qP}_R
 \le p!q!\norm P_{R+1}.
\]
Applying the tail estimates at radius $R+1$ constructs the derivative
codes uniformly. Taking $R=1$ gives the asserted derivative bound.
All limits used in these constructions have specified convergence rates.
\end{proof}

Exact rational coefficients will be obtained below by arranging that
each coefficient receives contributions from only finitely many stages.

\subsection{Orthogonal polynomials and tree flows}

\begin{lemma}[$\RCA$]\label{lem:positiveintegral}
If $g\in C([0,1])$ is nonnegative and $\Leb(g)=0$, then $g=0$.
\end{lemma}
\begin{proof}
If $g(a)>0$, a positive rational lower bound for $g(a)$ and the given
continuity modulus yield a nondegenerate rational subinterval on which
$g$ has a positive rational lower bound. This implies $\Leb(g)>0$.
\end{proof}

\begin{lemma}[$\RCA$]\label{lem:flowpath}
Let $T\subseteq2^{<\mathbb N}$ be prefix closed. If a jointly coded
family of reals $(a_\sigma)_\sigma$ satisfies
\[
 a_\varnothing=1,\qquad
 a_\sigma=a_{\sigma0}+a_{\sigma1},\qquad
 \sigma\notin T\Longrightarrow a_\sigma=0,
\]
then $T$ has a path. The values $a_\sigma$ need not be nonnegative.
\end{lemma}
\begin{proof}
Start at the root with $\eta=1/2$. If $|a_\sigma|>\eta$, approximate
both child values within $\eta/8$. Their sum is $a_\sigma$, so one
approximation has absolute value greater than $3\eta/8$. Choose that
child and replace $\eta$ by $\eta/4$. Its actual absolute value exceeds
the new lower bound; hence it belongs to $T$. Recursion using these
bounded rational comparisons gives a path.
\end{proof}

\begin{lemma}[$\RCA$]\label{lem:skewbasis}
There is a uniformly constructible sequence of nonzero rational
polynomials $V_n$ such that
\[
 V_n\in\operatorname{span}_{\mathbb Q}\{x^{n+1},\ldots,x^{2n+2}\},
 \qquad |V_n|_{\mathrm c}=1,
\]
and
\[
 \Leb(V_n)=0,\qquad \Leb(V_nV_k)=0\quad(k<n).
\]
Consequently $\zeta_n=\Leb(V_n^2)$ is a positive rational number,
and $V_n$ vanishes to order at least $n+1$ at the origin.
\end{lemma}
\begin{proof}
After constructing $V_0,\ldots,V_{n-1}$, impose orthogonality to
$1,V_0,\ldots,V_{n-1}$ on the rational space spanned by
$x^{n+1},\ldots,x^{2n+2}$. These are $n+1$ homogeneous rational
linear equations in $n+2$ variables. Rational Gaussian elimination
produces a nonzero solution, which we divide by its coefficient norm.
This is a recursion on finite rational codes. The square integral is
rational and positive, since a nonzero rational polynomial is nonzero
at a rational point and has an explicit continuity bound.
\end{proof}

\begin{lemma}[$\RCA$]\label{lem:skewdiagonal}
Every mixed equilibrium $(\mu,\nu)$ of a skew-symmetric zero-sum game
has value zero. Moreover, $(\mu,\mu)$ and $(\nu,\nu)$ are equilibria,
and these profiles are obtained uniformly from $(\mu,\nu)$.
\end{lemma}
\begin{proof}
Skew-symmetry gives $A(\rho,\rho)=0$ for every probability functional
$\rho$. If $(\mu,\nu)$ has value $v$, a row deviation to $\nu$ gives
$v\ge0$, and a column deviation to $\mu$ gives $v\le0$. Thus $v=0$.
The row inequalities give $A(x,\nu)\le0$ for all $x$, and skew-symmetry
then gives $A(\nu,y)\ge0$ for all $y$. Hence $(\nu,\nu)$ is an
equilibrium. The column inequalities give the same conclusion for
$(\mu,\mu)$.
\end{proof}

\subsection{The reversal}

\begin{proof}[Proof of Theorem~\ref{thm:skewequilibrium}]
The forward implication follows from Theorem~\ref{thm:main} and
Lemma~\ref{lem:skewdiagonal}. For the converse, let $T$ be an infinite
prefix-closed binary tree. Using Lemma~\ref{lem:skewbasis}, put
\[
 E_n=V_{2n},\qquad O_j=V_{2j+1},\qquad D_n=4n+2,
 \qquad b_n=2^{-n-(D_n+4)^2},\qquad
 \vartheta_n=b_n\zeta_{2n}.
\]
The two families are mutually orthogonal and consist of mean-zero
polynomials of coefficient norm one. Enumerate binary strings by
$\sigma_0=\varnothing$, $\sigma_{2n+1}=\sigma_n0$, and
$\sigma_{2n+2}=\sigma_n1$. Define
\begin{align}
 r_0&=E_0-\vartheta_0,\nonumber\\
 r_{2n+1}&=\frac{E_n}{\vartheta_n}
          -\frac{E_{2n+1}}{\vartheta_{2n+1}}
          -\frac{E_{2n+2}}{\vartheta_{2n+2}},\label{eq:skewfeatures}\\
 r_{2n+2}&=\begin{cases}
 E_n/\vartheta_n,&\sigma_n\notin T,\\
 0,&\sigma_n\in T.
 \end{cases}\nonumber
\end{align}
Let $f_j=\Nrm(r_j)$, and choose a degree bound $d_j$ for both $f_j$
and $O_j$. Set
\begin{equation}\label{eq:skewkernel}
 w_j=2^{-j-(d_j+4)^2},\qquad
 A_T(x,y)=\sum_{j\ge0}w_j
 \bigl(f_j(x)O_j(y)-O_j(x)f_j(y)\bigr).
\end{equation}
The payoff is skew-symmetric, and Lemma~\ref{lem:entire} gives
\begin{equation}\label{eq:skewtail}
 \norm{A_T-A_{T,<s}}_R\le2^{R^2-s+2},\qquad
 \norm{A_T}_2\le2^{-14}.
\end{equation}
In particular, the common bounds defining $\Aclass$ hold.

Since $O_j$ vanishes to order at least $2j+2$, the coefficient of
$x^ay^b$ in \eqref{eq:skewkernel} receives contributions only when
\[
 2j+2\le\max\{a,b\}.
\]
Every coefficient is therefore a uniformly computable finite rational
sum. This also applies to the constant term in $r_0$, which is paired
with $O_0$. For an integral radius $R\ge1$ and precision $2^{-k}$,
take
\[
 s=R^2+k+3,\qquad M=2\max_{j<s}d_j.
\]
All polynomial stages $j<s$ have total degree at most $M$. The norm
at radius $R$ of the Taylor terms of total degree greater than $M$
is therefore at most $2^{-k-1}$, by \eqref{eq:skewtail}. Thus the
exact coefficient array has a joint Taylor-truncation modulus, and
$A_T\in\Aclass^{\mathbb Q}$.

For each $s$, choose the lexicographically first node $\tau_s\in T$
of length $s+1$, and put
\[
 h_s=1+\sum_{\sigma_n\preceq\tau_s}b_nE_n,
 \qquad \rho_s(g)=\Leb(gh_s).
\]
Since $\sum_nb_n<1/4$, these are uniformly coded rational polynomial
densities with
\[
 \frac34\le h_s\le\frac54,\qquad \Leb(h_s)=1.
\]
Orthogonality gives $\rho_s(O_j)=0$ for every $j$, while
$\rho_s(E_n)/\vartheta_n$ is the indicator of the prefixes of $\tau_s$.
For a flow constraint $j=2n+1<s$, we have $|\sigma_n|\le n<s+1$,
so its parent lies strictly above the terminal level of $\tau_s$.
Excluded nodes have zero moment by prefix closure.
Hence $\rho_s(f_j)=0$ for $j<s$, including the root constraint.
The radius-two tail estimate in Lemma~\ref{lem:entire} now gives
\[
 \sup_y|A_T(\rho_s,y)|\le2^{-s-14}.
\]

Let $(\mu,\mu)$ be a symmetric equilibrium and put
$u(x)=A_T(x,\mu)$. Its value is zero, so $u\le0$. Therefore
\[
 0\le\frac34\Leb(-u)
 \le\rho_s(-u)
 =-A_T(\rho_s,\mu)
 \le2^{-s-14}.
\]
It follows that $\Leb(-u)=0$, and Lemma~\ref{lem:positiveintegral}
yields $u=0$.

Each $f_k$ belongs to the span of $1,E_0,E_1,\ldots$, whence
$\Leb(O_jf_k)=0$. Multiplying the zero section by $O_j$ and
integrating, with the uniform tail justifying passage to the series,
gives
\begin{equation}\label{eq:skewextract}
 0=\Leb\bigl(O_j(x)A_T(x,\mu)\bigr)
   =-w_j\zeta_{2j+1}\mu(f_j).
\end{equation}
Since $w_j\zeta_{2j+1}>0$, every constraint moment vanishes. Thus
\[
 a_n=\frac{\mu(E_n)}{\vartheta_n}
\]
satisfies
\[
 a_0=1,\qquad
 a_n=a_{2n+1}+a_{2n+2},\qquad
 \sigma_n\notin T\Longrightarrow a_n=0.
\]
Lemma~\ref{lem:flowpath} gives a path through $T$, proving $\WKL$.
No nonnegativity of these flow values is required.

Conversely, a path $P$ determines the density
\[
 h_P=1+\sum_{\sigma_n\sqsubset P}b_nE_n,
\]
which has integral one and lies in $[3/4,5/4]$. For each integer
$R\ge1$,
\[
 \sum_{n\ge N}b_n\norm{E_n}_R\le2^{R^2-N+1}.
\]
Thus $h_P$ has an effective entire-function code, and its coefficients
are rational by the increasing orders of vanishing. Its probability
functional annihilates every $f_j$ and $O_j$, so its diagonal profile
is an equilibrium. Finally, Lemma~\ref{lem:skewdiagonal} allows an
arbitrary equilibrium to replace a symmetric one. Since the constructed
payoffs have rational coefficients, the same reversal also applies to
$\Aclass$.
\end{proof}

\begin{remark}
Equation~\eqref{eq:skewextract} determines the constraint moments but
need not determine all $O_j$-moments, since the $f_j$ may be linearly
dependent. No description of the full equilibrium set as a common
annihilator is used. The infinite-tree hypothesis supplies the positive
approximate equalizers $\rho_s$.
\end{remark}

\begin{corollary}\label{cor:seqskew}
Over $\RCA$, $\WKL$ is equivalent to the assertion that every
uniformly coded sequence of games with payoffs in $\Aclass^{\mathbb Q}$
has a sequence of symmetric mixed Nash equilibria.
\end{corollary}
\begin{proof}
Apply Proposition~\ref{prop:upper} and the uniform diagonalization of
Lemma~\ref{lem:skewdiagonal}. The converse follows by constant
repetition and Theorem~\ref{thm:skewequilibrium}.
\end{proof}

\section{Countably many players}
\label{sec:countably-many-players}

The principle $\mathsf{SeqCNE}$ concerns a sequence of finite-player
games. In this section we consider a single game with countably many
players, whose payoffs may depend on the entire action profile.

\begin{definition}\label{def:countable-game}
A \emph{countable continuous game} consists of a uniformly represented
sequence of nonempty compact complete separable metric spaces
$(\widehat A_i,d_i)_{i\in\mathbb N}$ and a uniformly coded sequence
$f_i\in C(A)$, where
\[
 A=\prod_{i\in\mathbb N}\widehat A_i,
 \qquad
 D(x,y)=\sum_{i=0}^{\infty}2^{-i-1}
                  \min\{1,d_i(x_i,y_i)\}.
\]
The compact presentations are supplied jointly, as in
Definition~\ref{def:gamecode}, and the payoff codes include a jointly
coded sequence of moduli of uniform continuity with respect to $D$.
A mixed strategy profile is a jointly coded sequence
$(\mu_i)_{i\in\mathbb N}$ with $\mu_i\in\Prob(\widehat A_i)$.
\end{definition}

\begin{lemma}[$\RCA$]\label{lem:countable-products}
For a countable continuous game, the following constructions are uniform.
\begin{enumerate}
\item The space $(A,D)$ has a compact presentation and a specified
default profile $b=(b_i)_i$, where each $b_i$ belongs to the given dense
subset of $\widehat A_i$.
\item Every coded $f\in C(A)$ has a fast uniform approximation by
rational finite-cylinder tensor sums.
\item Every mixed strategy profile determines a positive unital
functional $\bigotimes_{i\in\mathbb N}\mu_i$ on $C(A)$ satisfying
\[
 \left(\bigotimes_i\mu_i\right)
       \left(\prod_{i\in J}p_i(x_i)\right)
       =\prod_{i\in J}\mu_i(p_i)
\]
for every finite $J\subseteq\mathbb N$ and every choice of basic
functions $p_i$.
\end{enumerate}
\end{lemma}

\begin{proof}
Choose $b_i$ uniformly from the given dense subsets. The profiles with
dense entries that agree with $b$ outside a finite initial segment form
a countable dense subset of $A$. For a $2^{-k}$-net, take
$2^{-k-2}$-nets in the first $N=k+2$ factors and fix the remaining
coordinates at $b$. The prefix error is at most $2^{-k-2}$, and the
tail error is at most $2^{-N}$. A fast $D$-Cauchy name gives a fast
Cauchy name in each coordinate by a computable shift of precision.
Taking the coordinate limits gives the limit in $D$, with a prescribed
rate. Thus the usual construction of products of coded metric spaces
\cite{simpson} is uniform in the given data.

For $x\in A$, put
\[
 x^{[N]}=(x_0,\ldots,x_{N-1},b_N,b_{N+1},\ldots).
\]
Since $D(x,x^{[N]})\le2^{-N}$, the supplied modulus of $f$ determines
$N$ such that $\norm{f-f(\cdot^{[N]})}_\infty<\eps/2$.
The finite-product approximation obtained from
\cite{ppy} gives a
rational tensor sum within $\eps/2$ of the restriction to the
first $N$ coordinates. This proves the second assertion.

Define the functional on finite-cylinder tensor sums by the
finite-product construction of \cite{ppy}.
Adding unit factors does not change its value. Moreover, a cylinder
function has the same uniform norm on $A$ as on its finite coordinate
domain. The finite-product contraction estimate therefore gives
\[
 \left|\left(\bigotimes_i\mu_i\right)(P)\right|
 \le\norm P_\infty.
\]
Consequently the definition is independent of the expression for $P$
and extends, using the approximations above, to a unital contraction on
$C(A)$. Positivity follows as in Lemma~\ref{lem:coordinates}.
The construction uses only finitely many marginal coordinates at each
prescribed precision.
\end{proof}

Write $E_\mu(f)=(\bigotimes_i\mu_i)(f)$, and write
$E_{a,\mu_{-i}}(f)$ for the value obtained by replacing $\mu_i$ with
$\delta_a$. Let $\mathsf{CNE}_\omega$ denote the assertion that every
countable continuous game has a mixed strategy profile satisfying
\[
 E_\mu(f_i)\ge E_{a,\mu_{-i}}(f_i)
 \qquad(i\in\mathbb N,\ a\in\widehat A_i).
\]

\begin{theorem}[$\RCA$]\label{thm:countable-player-equilibria}
The following are equivalent.
\begin{enumerate}
\item $\WKL$.
\item $\mathsf{CNE}_\omega$.
\item Every countable game with $A_i=[0,1]$ and payoffs
$f_i(x)=r_ix_i$, where $(r_i)_i$ is a jointly coded sequence in $[-1,1]$,
has a mixed Nash equilibrium.
\end{enumerate}
The equivalence remains valid if a pure equilibrium is required in the
third assertion.
\end{theorem}

\begin{proof}
Assume $\WKL$. Apply Lemma~\ref{lem:coordinates} simultaneously to all
players. In the resulting compact coordinate space, impose the
unitality, rational linearity, and contraction conditions for each
marginal. By Lemma~\ref{lem:countable-products}, choose uniformly
finite-cylinder tensor sums $P_{i,k}$ such that
$\norm{f_i-P_{i,k}}_\infty\le\eps_k=2^{-k}$. For every action $a$ in
the given dense subset of $\widehat A_i$, also impose
\[
 E_v(P_{i,k})\ge E_{a,v_{-i}}(P_{i,k})-2\eps_k.
\]
These are effective closed conditions, since the evaluations of the
tensor sums are polynomials in finitely many marginal
coordinates.

Consider a finite subfamily of these conditions. Let $I$ be the finite
set of players whose Nash inequalities occur. For each $i\in I$, take
a finite action list containing $b_i$ and all deviations specified in
the subfamily. Fix all players outside $I$ at their default actions.
The finite Nash theorem in $\WKL$, as observed in
\cite{ppy}, gives an equilibrium of this
explicitly listed finite game. Its atomic marginals, together with
$\delta_{b_j}$ for $j\notin I$, form a jointly coded profile. Its
expectations agree with those in the finite game, first on cylinder
tensor sums and then on $C(A)$ by uniform approximation. The selected
Nash conditions follow from the finite-game inequalities and the two
approximation errors. All functional conditions also hold. If $I$ is
empty, use the default Dirac profile. The closed conditions are
therefore finitely satisfiable.

The finite-intersection form of compactness in $\WKL$
\cite{simpson} gives a common solution. For every player $i$
and every dense action $a$, it satisfies
\[
 E_\mu(f_i)-E_{a,\mu_{-i}}(f_i)\ge-4\eps_k
 \qquad(k\in\mathbb N).
\]
Letting $k$ tend to infinity gives the exact inequality. The supplied
modulus extends it to every action. This proves the second assertion,
which implies the third.

For the reversal, let $\alpha,\beta:\mathbb N\to\mathbb N$ be
injections with disjoint ranges. Put
\[
 r_i=\sum_{t=0}^{\infty}2^{-t-2}
       \bigl(\mathbf1_{\alpha(t)=i}-\mathbf1_{\beta(t)=i}\bigr).
\]
The geometric tail gives uniform real codes. The function
$\omega_i(k)=i+k+2$ is a uniform continuity modulus for $r_ix_i$
with respect to $D$. Let $(\mu_i)_i$ be an equilibrium and set
$m_i=\mu_i(x\mapsto x)$. Deviation to $1$ gives $m_i=1$ whenever
$i\in\ran(\alpha)$, and deviation to $0$ gives $m_i=0$ whenever
$i\in\ran(\beta)$. Obtain uniformly rational approximations $q_i$
with $|q_i-m_i|<1/8$. By recursive comprehension,
$\{i:q_i>1/2\}$ is a set separating the two ranges. Hence $\WKL$
follows from $\Sigma^0_1$ separation \cite{simpson}.

For the pure version of the third assertion, $\WKL$ separates the
disjoint $\Sigma^0_1$ predicates $r_i>0$ and $r_i<0$. Choose action
$1$ on a separator and action $0$ elsewhere. The resulting profile is
a pure equilibrium. Conversely, Dirac evaluation turns a pure
equilibrium into a mixed equilibrium, so the same reversal applies.
\end{proof}

\begin{remark}
Replacing $r_ix_i$ by $2^{-i-1}r_ix_i$ preserves the best responses and
gives all payoffs the common Lipschitz bound $1$ with respect to $D$.
The use of real coefficient codes is essential to this reversal. For
exact rational coefficients, their signs are decidable and $\RCA$
constructs maximizing endpoints. The same conclusion holds for real
coefficients satisfying $r_i\ne0$ for every $i$, since a uniform search
for their signs then terminates. Thus the reversal does not establish
the result under uniqueness of every player's optimum. Throughout the
theorem, continuity means joint continuity for the stated product
metric, with the specified moduli; the proof does not apply to separate
continuity alone.
\end{remark}

\subsection{Approximate equilibria}

A rational finite-support mixed strategy has rational probability
weights on a finite list of specified dense action points. A mixed
profile has regret at most $\eps$, or is an $\eps$-Nash equilibrium,
if every unilateral deviation increases the corresponding expected
payoff by at most $\eps$.

\begin{lemma}[$\RCA$]\label{lem:rational-approximate-equilibria}
For every explicitly listed finite rational game and every positive
rational $\eta$, there is a rational mixed profile with regret less
than $\eta$. Such a profile is obtainable by a uniform numerical search.
\end{lemma}

\begin{proof}
For a finite rational table $\Gamma$, a positive rational $\eta$, and
a candidate rational mixed profile $q$, the condition
\[
 \operatorname{RatProfile}(\Gamma,q)\ \wedge\
 \bigwedge_{i,a}
 \bigl(E_{a,q_{-i}}(\Gamma_i)-E_q(\Gamma_i)<\eta\bigr)
\]
is decidable. The finite Nash theorem in $\WKL$
\cite{ppy}, followed by rational
approximation in the finite product of simplices, proves that such a
numerical code exists for every $\Gamma$ and $\eta$. This assertion is
arithmetic. Harrington's $\Pi^1_1$-conservation theorem
\cite{simpson} therefore makes it a theorem of $\RCA$.
Enumeration of rational profiles gives a uniform search whose
termination is provable in $\RCA$.
\end{proof}

\begin{corollary}[$\RCA$]\label{cor:base-system-approximation}
The following hold uniformly under the specified presentations.
\begin{enumerate}
\item Every finite-player continuous game has a rational finite-support
$\eps$-Nash equilibrium for each positive rational $\eps$.
\item Every countable continuous game has a jointly coded sequence of
profiles $(\mu^{(s)})_{s\in\mathbb N}$, each with rational
finite-support marginals, such that $\mu^{(s)}_j=\delta_{b_j}$ for
$j\ge s$ and
\[
 E_{a,\mu^{(s)}_{-i}}(f_i)-E_{\mu^{(s)}}(f_i)\le2^{-s}
 \qquad(i<s,\ a\in\widehat A_i).
\]
\item If the payoffs of a countable continuous game share a supplied
modulus of uniform continuity for $D$, then the game has a rational
finite-support $\eps$-Nash equilibrium for every positive rational
$\eps$, with only finitely many marginals different from the default
Dirac marginals. In particular, this holds under a supplied common
Lipschitz bound.
\end{enumerate}
In the first and third assertions, the profiles can be constructed
jointly for every prescribed sequence of positive rational tolerances.
\end{corollary}

\begin{proof}
For a finite-player game, choose finite nets of dense actions so fine
that replacing a player's action by a nearby net action changes that
player's payoff by at most $\eps/4$, uniformly in the other actions.
Approximate each payoff entry on the product grid by a rational number
within $\eps/8$. Lemma~\ref{lem:rational-approximate-equilibria}
gives a rational mixed profile with regret less than $\eps/4$ in the
resulting rational game. For a grid deviation, the two table errors add
at most $\eps/4$. Approximating an arbitrary deviation by a grid
action adds at most $\eps/4$. Thus the regret is less than $\eps$.
The construction is uniform, so recursive comprehension also gives the
profiles for a prescribed sequence of tolerances.

For the second assertion, let $\mu^{(0)}$ be the default Dirac
profile. At stage $s>0$, fix all coordinates $j\ge s$ at $b_j$ and
apply the first assertion to the game of the first $s$ players with
tolerance $2^{-s}$. By Lemma~\ref{lem:countable-products}, the finite
expectations agree with the original product-functional expectations.
Uniformity gives the required joint code.

For the third assertion, let $\omega$ be the common modulus. Choose
$k$ with $2^{-k}<\eps/2$ and put $N=\omega(k)+1$. A unilateral
change by player $j\ge N$ changes $D$ by at most $2^{-j-1}$ and hence
changes that player's payoff by less than $\eps/2$. Apply the first
assertion to the first $N$ players with tolerance $\eps/2$, fixing all
remaining coordinates at their default values. The resulting profile
has the required regret bound for the first $N$ players. For every
other player, the pointwise bound passes to expectations by positivity
of the product functional. This construction is again uniform in
$\eps$. Under a common Lipschitz bound $L$, it suffices to choose
$N$ with $L2^{-N-1}<\eps/2$.
\end{proof}

\begin{remark}
The common-modulus hypothesis cannot be omitted from the third
assertion. If $A_i=\{0,1\}$, $b_i=0$, and $f_i(x)=x_i$, then every
player left at the default action has regret one. On the other hand,
the rescaled affine games above have common Lipschitz bound $1$.
Within that class, prescribed-regret approximation is provable in
$\RCA$, whereas exact equilibrium existence is equivalent to $\WKL$.
The approximate profiles are not asserted to converge.
\end{remark}

\textbf{Declaration of generative AI use}

The author used ChatGPT (OpenAI) to assist with English-language editing, manuscript organization, and the review and clarification of mathematical arguments.

\end{document}